\documentclass[12pt, a4paper, twoside]{article}
\usepackage{amssymb}
\usepackage[leqno]{amsmath}
\usepackage{latexsym}
\usepackage{enumerate}
\usepackage{amsmath, amssymb, amscd}
\Roman{equation}

\DeclareMathOperator{\Pic}{Pic}
\DeclareMathOperator{\NS}{NS}

\DeclareMathOperator{\Supp}{Supp}

\begin{document}

\title{\textbf{\Large{On the classification of non-aCM curves on quintic hypersurfaces in $\mathbb{P}^3$}}}

\author{Kenta Watanabe \thanks{Nihon University, College of Science and Technology,   7-24-1 Narashinodai Funabashi city Chiba 274-8501 Japan , {\it E-mail address:watanabe.kenta@nihon-u.ac.jp}, Telephone numbers: 090-9777-1974} }

\date{}

\maketitle 

\begin{abstract}

\noindent In this paper, we call a sub-scheme of dimension one in $\mathbb{P}^3$ a curve. It is well known that the arithmetic genus and the degree of an aCM curve $D$ in $\mathbb{P}^3$ is computed by the $h$-vector of $D$. However, for a given curve $D$ in $\mathbb{P}^3$, the two invariants of $D$ do not tell us whether $D$ is aCM or not. In this paper, we give a classification of curves on a smooth quintic hypersurface in $\mathbb{P}^3$ which are not aCM.

\end{abstract}


\smallskip

\noindent {\bf{Mathematics Subject Classifcation }}14J29, 14J60, 14J70

\section{Introduction} 

We work over the complex number field $\mathbb{C}$. Let $C$ be a curve in $\mathbb{P}^n$ which is not necessarily irreducible. Then we call $C$ an  {\it{arithmetically Cohen-Macaulay}} ({{\it{aCM}} for short) curve if the homogeneous coordinate ring of it is Cohen-Macaulay. It is equivalent to the condition that $h^1(\mathcal{I}_C(l))=0$ for all $l\in\mathbb{Z}$, where $\mathcal{I}_C$ is the ideal sheaf of $C$ in $\mathbb{P}^n$. It seems that the projective normality of a smooth curve in $\mathbb{P}^n$ is often defined by the equivalent condition that the natural map $H^0(\mathcal{O}_{\mathbb{P}^n}(l))\longrightarrow H^0(\mathcal{O}_C(l))$ is surjective for any $l\geq0$. In general, it is difficult to discriminate whether a given curve $C$ in $\mathbb{P}^n$ is aCM or not. However, if $C$ lies on a smooth projective surface, we can do it in several cases. For examples, aCM curves on a DelPezzo surface $X$ with respect to the projective embedding induced by $-K_X$ are classified by Pons-Llopis and Tonini ([4]). On the other hand, it is well known that any non-hyperellipitic curve $C$ on a K3 surface has a canonical embedding and is projectively normal with respect to it (cf. [5]). A curve $C$ on a smooth hypersurface $X$ in $\mathbb{P}^3$ is aCM if and only if $\mathcal{O}_X(C)$ is an aCM line bundle on $X$.  In the previous work [6,7], aCM line bundles on a smooth hypersurface of degree $d=4,5$ in $\mathbb{P}^3$ are characterized. In particular, a necessary and sufficient condition for a line bundle on a smooth quintic hypersurface $X$ in $\mathbb{P}^3$ to be aCM is given as follows.

\newtheorem{thm}{Theorem}[section]

\begin{thm} {\rm{([7, Theorem 1.1])}} Let $X$ be a smooth quintic hypersurface in $\mathbb{P}^3$, let $H_X$ be the hyperplane class of $X$, and let $C$ be a smooth member of the linear system $|H_X|$. Let $D$ be a non-zero effective divisor on $X$ of arithmetic genus $P_a(D)$, and we set $k=C.D+1-P_a(D)$. Then $\mathcal{O}_X(D)$ is aCM and initialized if and only if the following conditions are satisfied.

$\;$

\noindent {\rm{(i)}} $0\leq k\leq 4$.

\smallskip

\smallskip

\noindent {\rm{(ii)}} If $0\leq k\leq 1$, then $C.D=10-k$ and $h^0(\mathcal{O}_C(D-C))=0$. 

\smallskip

\smallskip

\noindent {\rm{(iii)}} If $k=2$, then the following conditions are satisfied.

\smallskip

\smallskip

{\rm{(a)}} $C.D\in\{1,4,7,8\}$.

{\rm{(b)}} If $C.D=7$, then $h^0(\mathcal{O}_X(2C-D))=1$.

{\rm{(c)}} If $C.D=8$, then $h^0(\mathcal{O}_C(D-C))=0$ and $h^0(\mathcal{O}_C(D))=3$.

\smallskip

\smallskip

\noindent {\rm{(iv)}} If $3\leq k\leq 4$, then the following conditions are satisfied.

\smallskip

{\rm{(a)}} $k-1\leq C.D\leq 10-k$.

{\rm{(b)}} If $8-k\leq C.D\leq 10-k$, then $h^0(\mathcal{O}_C(D))=5-k$. 

{\rm{(c)}} If $k=3$, then $C.D\neq4$. 

\smallskip

\smallskip

\noindent Here, we say that a line bundle $L$ on $X$ is initialized if $|L|\neq\emptyset$ and $|L(-1)|=\emptyset$. \end{thm} 

The Hilbert function of an aCM curve $C$ in $\mathbb{P}^3$ is defined by $H_C(l)=h^0(\mathcal{O}_C(l))$ ($l\in\mathbb{N}\cup\{0\}$), and the $h$-vector of $C$ is defined as the second difference function of it. The degree $\deg C$ and the arithmetic genus $P_a(C)$ of $C$ are denoted as follows respectively.
$$\deg C=\displaystyle\sum_{l\geq0}h_C(l)\;\;\;\;\;\;\;\;\;\;\; P_a(C)=\displaystyle\sum_{l\geq1}(l-1)h_C(l).$$
However, we can not discriminate whether a given curve in $\mathbb{P}^n$ is aCM or not by using the two invariants of it. For example, any smooth irreducible aCM curve of degree 6 and genus 3 in $\mathbb{P}^3$ is not hyperelliptic and a smooth quartic hypersurface in $\mathbb{P}^3$ containing such an aCM curve is linear determinantal. However, a smooth curve of bidegree $(2,4)$ on a smooth quadric hypersurface in $\mathbb{P}^3$ is a hyperelliptic curve which has the same degree and genus as it. Hence, it is natural and interesting to classify the pair $(d,g)$ consisting of integers $d$ and $g$ such that there exist an aCM curve of degree $d$ and genus $g$ and  a non-aCM curve of the same degree and genus as it in $\mathbb{P}^3$. However, any concrete description concerning non-aCM curves on $X$ is not given in Theorem 1.1. In this paper, we will give a best possible sufficient condition regarding the degree and the arithmetic genus, for a curve on a smooth quintic hypersurface $X$ in $\mathbb{P}^3$ to be aCM, and a complete classification of non-aCM curves on $X$. Our first main theorem is the following.

\begin{thm} Let the notation be as in Theorem 1.1. If $D$ satisfies one of the following conditions, then $D$ is an aCM curve.

$\;$

\noindent{\rm{(i)}} $k=2$ and $C.D\in \{1,4\}$.

\noindent{\rm{(ii)}} $k=3$ and $C.D\in \{2,3,5,6\}$.

\noindent{\rm{(iii)}} $k=4$ and $C.D\in \{3,4\}$. \end{thm}

\noindent By the assertion of Theorem 1.1, the essence of Theorem 1.2 is that if $(k,C.D)=(3,5),(3,6)$ or $(4,4)$, then $D$ is aCM. On the other hand, since the existence of an aCM curve $D$ on a smooth quintic hypersurface in $\mathbb{P}^3$ satisfying each condition in Theorem 1.1 is indicated in Section 5 of [7], we can classify the pair consisting of integers $d$ and $g$ such that the Hilbert scheme of curves of degree $d$ and arithmetic genus $g$ in $\mathbb{P}^3$ contains both an aCM curve and a non-aCM curve on a smooth quintic hypersurface in $\mathbb{P}^3$. Our second main theorem is the following.

\begin{thm} Let $d$ and $k$ be integers. If $d$ and $k$ satisfy one of the following conditions, then there exists a non-aCM curve $D$ of degree $d$ and arithmetic genus $P_a(D)=d-k+1$ on a smooth quintic hypersurface in $\mathbb{P}^3$.

$\;$

\noindent{\rm{(i)}} $k=0$ and $d=10$.

\noindent{\rm{(ii)}} $k=1$ and $d=9$.

\noindent{\rm{(iii)}} $k=2$ and $d\in \{7,8\}$.

\noindent{\rm{(iv)}} $k=3$ and $d=7$.

\noindent{\rm{(v)}} $k=4$ and $d\in\{5,6\}$. \end{thm}

\noindent Non-aCM curves satisfying each condition in Theorem 1.3 can be concretely constructed, and we can obtain a complete classification of non-aCM curves on a smooth quintic hypersurface in $\mathbb{P}^3$ (see Section 4).

Our plan of this paper is the following. In section 2, we recall the previous work about the classification of aCM line bundles on a smooth quartic hypersurface in $\mathbb{P}^3$, and give a characterization and examples of non-aCM curves on such a surface. In section 3, we recall several facts concerning line bundles on smooth quintic hypersurfaces  in $\mathbb{P}^3$. In section 4, we prove the main results and give a concrete description of non-aCM curves on smooth quintic hypersurfaces in $\mathbb{P}^3$.

$\;$

{\bf{Notations and conventions}}. In this paper, a surface is smooth and projective. Let $X$ be a smooth curve or a surface. Then we denote the canonical bundle of $X$ by $K_X$. For a divisor or a line bundle $L$ on $X$, we denote by $|L|$ the linear system of $L$, and denote the dual of a line bundle $L$ by $L^{\vee}$. For an irreducible curve $D$ which is not necessarily smooth, we denote by $P_a(D)$ the arithmetic genus of $D$. For irreducible divisors $D_1$ and $D_2$ on a surface $X$, the arithmetic genus of $D_1+D_2$ is denoted as $P_a(D_1)+P_a(D_2)+D_1.D_2-1$. By induction, the arithmetic genus of a reducible effective divisor $D$ on $X$  is also defined. We denote it by the same notation $P_a(D)$. It follows that $2P_a(D)-2=D.(K_X+D)$ by the adjunction formula. Note that if $D$ is reduced and irreducible, then $P_a(D)\geq0$.

The gonality of a smooth curve is the minimal degree of pencils on it. It is well known that the gonality of a smooth plane curve of degree $d\geq5$ is $d-1$. 

We denote the N\'{e}ron-Severi lattice and the Picard lattice of a surface $X$ by $\NS(X)$ and $\Pic(X)$ respectively. We call the rank of the N\'{e}ron-Severi lattice of $X$ the Picard number of $X$. If the Picard number of $X$ is $\rho$, then by the Hodge index theorem, the signature of $\NS(X)$ is $(1,\rho-1)$. This implies that $D_1^2D_2^2\leq(D_1.D_2)^2$ for two divisors $D_1$ and $D_2$ on $X$ with $D_1^2>0$ and $D_2^2>0$. 

Let $X$ be a smooth hypersurface of degree $d$ in $\mathbb{P}^3$, and let $C$ be a hyperplane section of $X$.  For a non-zero effective divisor $D$ on $X$, we denote the degree of $D$ by $\deg D:=C.D$. In particular, we call a curve on $X$ of degree 1 and genus 0 a line. We denote the class of $C$ in $\Pic(X)$ by $H_X$. For an integer $l$, $H_X^{\otimes l}$ is often denoted as $\mathcal{O}_X(l)$. By the adjunction formula, $K_X\cong\mathcal{O}_X(d-4)$. For a line bundle $L$ on $X$, we will write $L\otimes\mathcal{O}_X(l)=L(l)$.

\section{ACM curves on quartic hypersurfaces in $\mathbb{P}^3$}

In this section, we recall the result concerning aCM curves on smooth quartic hypersurfaces in $\mathbb{P}^3$, and give a first example of a non-aCM curve as a toy model of our work. First of all, we recall the relationship between aCM curves and aCM line bundles on a hypersurface in $\mathbb{P}^3$.

\newtheorem{df}{Definition}[section]

\begin{df} Let $X$ be a hypersurface in $\mathbb{P}^3$, and let $L$ be a line bundle on $X$. Then we say that $L$ is {\rm{arithmetically Cohen-Macaulay (aCM}} for short{{\rm{)}}} if $h^1(L(l))=0$ for any $l\in\mathbb{Z}$. \end{df}

\noindent{\bf{Remark 2.1}} Let $D$ be a non-zero effective divisor on a hypersurface $X$ of degree $d$ in $\mathbb{P}^3$. Then, by the exact sequence
$$0\longrightarrow\mathcal{O}_{\mathbb{P}^3}(-d)\longrightarrow\mathcal{I}_{D}\longrightarrow\mathcal{O}_X(-D)\longrightarrow0,$$
$\mathcal{O}_X(D)$ is aCM if and only if $D$ is aCM, where $\mathcal{I}_{D}$ is the ideal sheaf of $D$ in $\mathbb{P}^3$.

$\;$

\noindent By Remark 2.1, our previous work as in [6] regarding the classification of aCM line bundles on a quartic hypersurface $X$ in $\mathbb{P}^3$ implies the following assertion.

\newtheorem{prop}{Proposition}[section]

\begin{prop} Let $X$ be a quartic hypersurface in $\mathbb{P}^3$, and let $D$ be a non-zero effective divisor on $X$. Let $C$ be a hyperplane section of $X$. Then the following conditions are equivalent.

\smallskip

\smallskip

\noindent {\rm{(i)}} $D$ is an aCM curve with $|D-C|=\emptyset$.

\noindent {\rm{(ii)}} One of the following cases occurs.

\smallskip

\smallskip

{\rm{(a)}} $P_a(D)=0$ and $1\leq C.D\leq 3$.

{\rm{(b)}} $P_a(D)=1$ and $3\leq C.D\leq 4$.

{\rm{(c)}} $P_a(D)=2$ and $C.D=5$.

{\rm{(d)}} $P_a(D)=3,\; C.D=6$, and $|D-C|=|2C-D|=\emptyset$. \end{prop}

\noindent It turns out that if a curve $D$ with $P_a(D)\leq 2$ belongs to a smooth quartic hypersurface in $\mathbb{P}^3$, we can discriminate whether or not it is aCM, by using the arithmetic genus and the degree of it. However, if a smooth curve $D$ with $\deg D=6$ and $P_a(D)=3$ in $\mathbb{P}^3$ is hyperelliptic, then it  is not aCM. In fact, if $D$ is aCM, then there exists a quartic hypersurface $X$ in $\mathbb{P}^3$ containing $D$ which is linear determinantal, and $|D|$ defines a birational map from $X$ onto the image of $X$. Since $X$ is a K3 surface, this implies that $K_D$ is very ample and hence, $D$ is not hyperelliptic. This means that if we remove the condition that $|D-C|=|2C-D|=\emptyset$ from Proposition 2.1 (ii) (d), the assertion is not correct.

\begin{prop} Let $X$, $C$ and $D$ be as in Proposition 2.1, and assume that $C.D=6$ and $P_a(D)=3$. Then the following conditions are equivalent.

\smallskip

\smallskip

\noindent {\rm{(a)}} $D$ is not an aCM curve.

\noindent {\rm{(b)}} There exist two lines $\Gamma_1$ and $\Gamma_2$ on $X$ with $\Gamma_1.\Gamma_2=0$ such that $\Gamma_1+\Gamma_2\in |D-C|$ or $\Gamma_1+\Gamma_2\in |2C-D|$. \end{prop}

\smallskip

\smallskip

{\it{Proof}}. (b) $\Longrightarrow$ (a) Let $\Gamma_1$ and $\Gamma_2$ be lines on $X$ with $\Gamma_1.\Gamma_2=0$. If $\Gamma_1+\Gamma_2\in |D-C|$, then $D-C$ is not 1-connected. Since $X$ is a K3 surface, we have $h^1(\mathcal{O}_X(D)(-1))=h^1(\mathcal{O}_X(1)(-D))\neq0$. This implies that $\mathcal{O}_X(D)$ is not aCM and hence, $D$ is not an aCM curve. If $\Gamma_1+\Gamma_2\in |2C-D|$, then $2C-D$ is not 1 connected. Hence, by the same reason as above, $D$ is not aCM. 

\smallskip

(a) $\Longrightarrow$ (b) We assume that $D$ is not aCM. By Proposition 2.1 and Remark 2.1, we have $|D-C|\neq\emptyset$ or $|2C-D|\neq\emptyset$. Assume that $|D-C|\neq\emptyset$. Since $C.(D-C)=2$ and $(D-C)^2=-4$, the arithmetic genus of the member of $|D-C|$ is $-1$ and hence, there exist two lines $\Gamma_1$ and $\Gamma_2$ on $X$ with $\Gamma_1.\Gamma_2=0$ such that $\Gamma_1+\Gamma_2\in |D-C|$. If $|2C-D|\neq\emptyset$, then there exist two lines $\Gamma_1$ and $\Gamma_2$ on $X$ with $\Gamma_1.\Gamma_2=0$ such that $\Gamma_1+\Gamma_2\in |2C-D|$, by the same reason as above. Hence, we have the assertion. $\hfill\square$

$\;$

\noindent In Proposition 2.2, if $|D-C|\neq\emptyset$, then the two lines $\Gamma_1$ and $\Gamma_2$ form the fixed component of $|D|$, and hence, we can not take such a divisor $D$ on $X$ to be smooth and  irreducible. On the other hand, if $|2C-D|\neq\emptyset$, then $D$ is linearly equivalent to the union of two elliptic curves $E_1$ and $E_2$ of degree 3 on $X$ with $E_1.E_2=2$ and hence, we can take such a curve $D$ to be smooth and irreducible. Since the restrictions to $D$ of the pencils $|E_1|$ and $|E_2|$ on $X$ are pencils of degree 2 on $D$, the curve $D$ is hyperelliptic. We can construct a non-aCM curve lying on the Fermat quartic hypersurface in $\mathbb{P}^3$ satisfying the condition (b).

$\;$

\noindent{\bf{Example 2.1}} Let $X$ be the quartic hypersurface in $\mathbb{P}^3$ defined by the equation $x_0^4+x_1^4+x_2^4+x_3^4=0$, for a suitable homogeneous coordinate $(x_0:x_1:x_2:x_3)$ on $\mathbb{P}^3$. Let $C_i$ be the hyperplane section of $X$ defined by the equation $x_0+\omega x_i=0$ ($i=1,2$), where $\omega$ is a primitive eighth root of unity. Let $\Gamma_1$ and $\Gamma_2$  be lines defined by the equations $x_0+\omega x_1=x_2+\omega x_3=0$ and $x_0+\omega x_2=x_1+\omega^2 x_3=0$, respectively. Since $\Gamma_1.\Gamma_2=0$, if we set $D_1=C_1+\Gamma_1+\Gamma_2$ and $D_2=C_1+C_2-\Gamma_1-\Gamma_2$, then $P_a(D_i)=3$ and $\deg D_i=6$ for $i=1,2$. Moreover, the conditions that $|D_1-C_1|\neq\emptyset$ and $|2C_1-D_2|\neq\emptyset$ are also satisfied.

\section{Linear systems on quintic hypersurfaces in $\mathbb{P}^3$}

From now on let $X$ be a smooth quintic hypersurface in $\mathbb{P}^3$. In this section, we recall several basic facts regarding linear systems on $X$ and some useful propositions as in Section 2 of [7]. Let $D$ be a divisor on $X$ and let $C$ be a smooth hyperplane section of $X$. First of all, since $K_X\cong\mathcal{O}_X(1)$, the Riemann-Roch theorem for $\mathcal{O}_X(D)$ implies that
$$\chi(\mathcal{O}_X(D))=\displaystyle\frac{1}{2}D.(D-C)+5,$$
where $\chi(\mathcal{O}_X(D))=h^0(\mathcal{O}_X(D))-h^1(\mathcal{O}_X(D))+h^2(\mathcal{O}_X(D))$. Note that since $h^0(K_X)=h^0(\mathcal{O}_{\mathbb{P}^3}(1))=4$, $h^1(\mathcal{O}_X)=0$, and $h^0(\mathcal{O}_X)=1$, we have $\chi(\mathcal{O}_X)=5$.

\noindent The Serre duality for $\mathcal{O}_X(D)$ is given by 
$$h^i(\mathcal{O}_X(D))=h^{2-i}(\mathcal{O}_X(1)(-D))\;(0\leq i\leq 2).$$

\noindent{\bf{Remark 3.1}} Since $C$ is a plane quintic curve, the gonality of $C$ is 4. Hence, a line bundle $L$ on $C$ with $h^0(L)\geq2$ satisfies $\deg(L)\geq h^0(L)+2$. In particular, if $|D-C|=\emptyset$ and $h^0(\mathcal{O}_X(D))\geq2$, then, by the exact sequence
$$0\longrightarrow\mathcal{O}_X(D)(-1)\longrightarrow\mathcal{O}_X(D)\longrightarrow\mathcal{O}_C(D)\longrightarrow0,\leqno (3.1)$$
we have $h^0(\mathcal{O}_X(D))\leq h^0(\mathcal{O}_C(D))\leq C.D-2$.

$\;$

\noindent By the Riemann-Roch theorem, the Serre duality, and Remark 3.1, a non-zero effective divisor on $X$ of degree $d\leq2$ can be explicitly written as follows.

\newtheorem{lem}{Lemma}[section]

\begin{lem} {\rm{([7, Lemma 2.1])}} Let $D$ be a divisor on $X$ satisfying $C.D=1$. Then the following conditions are equivalent.

\smallskip

\smallskip

\noindent{\rm{(a)}} $h^0(\mathcal{O}_X(D))>0$.

\smallskip

\noindent{\rm{(b)}} $h^0(\mathcal{O}_X(D))=1$, $h^0(\mathcal{O}_X(C-D))=2$, and $h^1(\mathcal{O}_X(D))=0$.

\smallskip

\noindent{\rm{(c)}} $D^2=-3$.\end{lem}

\begin{lem} {\rm{([7, Lemma 2.2])}} Let $D$ be an effective divisor on $X$ with $C.D=2$. If $D^2\leq-6$, then one of the following cases occurs.

\smallskip

\noindent{\rm(a)} There exists a line $D_1$ on $X$ with $D=2D_1$.

\smallskip

\noindent{\rm(b)} There exist two lines $D_1$ and $D_2$ such that $D=D_1+D_2$ and $D_1.D_2=0$.\end{lem}

\noindent The arithmetic genus $P_a(D)$ of a non-zero effective divisor $D$ is given as follows.
$$P_a(D)=\displaystyle\frac{1}{2}D.(D+C)+1.$$

\noindent{\bf{Remark 3.2}} Obviously, the arithmetic genus $P_a(D)$ of an effective divisor $D$ of degree one on $X$ is 0. On the other hand,  if a divisor $D$ on $X$ satisfies $C.D=2$, then by Remark 3.1, $h^0(\mathcal{O}_X(D))\leq1$ and $h^0(\mathcal{O}_X(1)(-D))\leq1$. Hence, by the Riemann-Roch theorem, we have $D^2\leq-4$. If $D^2=-4$, then the member of $|D|$ is a plane conic and hence, the arithmetic genus of it is zero.

$\;$

\noindent By Lemma 3.1 and Lemma 3.2, we have the following assertion {\rm{(cf. [7, Proposition 2.1])}}.

\begin{prop} Let $k$ be an integer with $0\leq k\leq 4$, and let $D$ be a non-zero effective divisor on $X$ such that $P_a(D)=C.D+1-k$. If $C.D\geq 7-k$, then $h^0(\mathcal{O}_X(1)(-D))=0$. \end{prop}

\noindent{\bf{Remark 3.3}} By the Riemann-Roch theorem and the Serre duality, if $D$ is a non-zero effective divisor on $X$ with $h^1(\mathcal{O}_X(-D))=0$, then $P_a(D)\geq0$.

$\;$

\noindent In general, the vanishing condition of the cohomology of the sheaf as in Remark 3.3 can be characterized by the following notion.

\begin{df} Let $m$ be a positive integer. Then a non-zero effective divisor $D$ on $X$ is called {\rm{$m$-connected}} if $D_1.D_2\geq m$, for each non-trivial effective decomposition $D=D_1+D_2$.\end{df}
\noindent A non-zero effective divisor $D$ is 1-connected if and only if $h^0(\mathcal{O}_D)=1$. Moreover, by the exact sequence
$$0\longrightarrow \mathcal{O}_X(-D)\longrightarrow \mathcal{O}_X\longrightarrow\mathcal{O}_D\longrightarrow 0,$$
the condition that $h^0(\mathcal{O}_D)=1$ is equivalent to $h^1(\mathcal{O}_X(-D))=0$. Hence, if $D$ is 1-connected, then we have $P_a(D)\geq0$, by Remark 3.3. Conversely, the following fact is well known as a sufficient condition for a non-zero effective divisor $D$ on $X$ with $P_a(D)\geq0$ to be 1-connected.

\begin{prop} {\rm{([1, p 179, Theorem 12.1])}} Let $D$ be a numerical effective divisor on $X$ with $D^2>0$. Then $h^1(\mathcal{O}_X(-D))=0$. \end{prop}

The effective divisor $D$ consisting of two skew lines on $X$ as in Lemma 3.2 (b) is not 1-connected, and is not contained in any hyperplane of $\mathbb{P}^3$. Since such a divisor $D$ satisfies $|D-C|=\emptyset$, by the Riemann-Roch theorem and Remark 3.1, the conditions $h^0(\mathcal{O}_X(D))=1$ and $h^1(\mathcal{O}_X(D))=0$ are also satisfied. Conversely, any non-zero effective divisor $D$ which is not 1-connected is characterized as follows, under the condition that $|D-C|=\emptyset$  and $h^1(\mathcal{O}_X(D))=0$ {\rm{(cf. [7, Proposition 2.3])}}. 

\begin{prop} Let $D$ be a non-zero effective divisor on $X$. If $|D-C|=\emptyset$ and $h^1(\mathcal{O}_X(D))=0$, then $h^1(\mathcal{O}_X(-D))=0$ or there exist two lines $D_1$ and $D_2$ on $X$ such that $D=D_1+D_2$ and $D_1.D_2=0$.\end{prop}

\noindent We can construct an example of an effective divisor $D$ on a smooth quintic hypersurface $X$ in $\mathbb{P}^3$ satisfying the condition (b) as in Lemma 3.2.

$\;$

\noindent {\bf{Example 3.1}} Let $X$ be the quintic hypersurface in $\mathbb{P}^3$ defined by the equation $x_0^5+x_1^5+x_2^5+x_3^5=0$. Let $D_1$ and $D_2$ be divisors on $X$ which are defined by the equations $x_0+x_1=x_2+x_3=0$ and $x_0+x_2=x_1+\xi x_3=0$, respectively. Here, $\xi$ is a primitive fifth root of unity. Then $D_1$ and $D_2$ are skew lines on $X$. Hence, the divisor $D_1+D_2$ on $X$ is not 1-connected.

$\;$

\noindent By Proposition 3.3, we have the following assertion as a sufficient condition for a non-zero effective divisor $D$ on $X$ to be 1-connected.

$\;$

\noindent{\bf{Corollary 3.1}} Let $D$ be a non-zero effective divisor on $X$ satisfying the condition that  $h^1(\mathcal{O}_X(D))=0$ and $|D-C|=\emptyset$. If $P_a(D)\geq0$, then $h^1(\mathcal{O}_X(-D))=0$.

$\;$

Assume that the linear system $|D|$ defined by a non-zero effective divisor $D$ is base point free and $\dim|D|=r$. Let $f:X\longrightarrow \mathbb{P}^r$ be the map defined by $|D|$. The theorem of Bertini implies that if the dimension of the image of $f$ is two, then the general member of $|D|$ is smooth and irreducible; otherwise, there exists a Stein factorization $f=h\circ g$ which is the composition of a pencil $g:X\longrightarrow \Gamma$ and a finite map $h:\Gamma\longrightarrow\mathbb{P}^r$ onto the image, where $\Gamma$ is a smooth irreducible curve. Here, we have $\Gamma\cong\mathbb{P}^1$. Indeed, since $h^1(\mathcal{O}_X)=0$, $\NS(X)\cong\Pic(X)$ and hence, the genus of the image of $g$ is 0. Hence, in the latter case, the general member of $|D|$ is the disjoint union of finitely many smooth curves $D_1,\cdots,D_m$ such that $D_i$ and $D_j$ are linearly equivalent for each $i$ and $j$. If $D$ is 1-connected, then this case does not occur. By Proposition 3.2, we have the following assertion.

\begin{prop} Let $D$ be a base point free divisor on $X$ with $D^2>0$. Then the general member of $|D|$ is a smooth irreducible curve. \end{prop}

\noindent{\bf{Remark 3.4}} For a non-zero effective divisor $D$ on $X$, if $|D|$ has no fixed component and $D^2=0$, then there exist a smooth irreducible curve $\tilde{D}$ and $m\in\mathbb{N}$ with $m\tilde{D}\in|D|$.

\smallskip

\smallskip

\noindent If $|D|$ has no fixed component, then the general member of it is smooth outside the set of base points of $|D|$. In particular, the following assertion follows.

\begin{prop} Let $D$ be an effective divisor on $X$ with $0<D^2\leq3$ such that $|D|$ has no fixed component. Then $|D|$ contains a smooth irreducible curve. \end{prop}

{\it{Proof}}. Assume that $|D|$ has a base point $P$ such that any member of $|D|$ is singular at $P$. Let $m$ be the minimum number of the multiplicity at $P$ of all members of $|D|$. Then $m\geq2$. Let $\varphi:\tilde{X}\longrightarrow X$ be a blow up at $P$ and we set $E=\varphi^{-1}(P)$. If we set $\tilde{D}=\varphi^{\ast}D-mE$, then $\tilde{D}^2\leq3-m^2\leq-1$. However, this contradicts the fact that $|\tilde{D}|$ has no fixed component. By the assumption, $D$ is not zero, and the cardinality of the set of the base points of $|D|$ is at most 3. Hence, the general member of $|D|$ is a smooth curve. Since $D$ is numerical effective and $D^2>0$, by Proposition 3.2, $h^1(\mathcal{O}_X(-D))=0$. Hence, the general member of $|D|$ is irreducible. $\hfill\square$

\section{Proof and examples of main Theorems}

In this section, we prove Theorem 1.2 and Theorem 1.3. Let $X$ be a quintic hypersurface in $\mathbb{P}^3$ and let $D$ be a non-zero effective divisor on $X$. Let $C$ be a smooth hyperplane section of $X$. We recall the following assertion which follows from Proposition 3.2 in [7] and Remark 2.1 to prove our main theorems.

\begin{prop} Let $k$ be a positive integer with $C.D+5<5k$. If $h^1(\mathcal{O}_X(l)(-D))=0$ for $0\leq l\leq k$, then $D$ is an aCM curve. \end{prop}

\noindent  First of all, we investigate the case where $D$ is an aCM curve. By Theorem 1.1, it is sufficient to consider the case where $(C.D,P_a(D))=(4,1),(5,3), and (6,4)$.

\begin{prop} If $P_a(D)=1$ and $C.D=4$, then $D$ is aCM. \end{prop}

{\it{Proof}}. Since $(D-C).C=-1$, we have $|D-C|=\emptyset$. By Proposition 4.1, it is sufficient to show that $h^1(\mathcal{O}_X(l)(-D))=0$, for $0\leq l\leq2$. We show that $h^1(\mathcal{O}_X(1)(-D))=0$. First of all, by Proposition 3.1, we have $|C-D|=\emptyset$. Since $\chi(\mathcal{O}_X(D))=1$, if $h^1(\mathcal{O}_X(1)(-D))\neq0$, then $h^0(\mathcal{O}_X(D))\geq2$. Let $\Delta$ be the fixed component of $|D|$ and we set $\tilde{D}=D-\Delta$. Since $|D-C|=\emptyset$, we have $|\tilde{D}-C|=\emptyset$. Since $D^2=-4<0$, $\Delta$ is not empty. By the ampleness of $C$, we have $C.\tilde{D}\leq3$. This contradicts Remark 3.1. Therefore, we have $h^1(\mathcal{O}_X(1)(-D))=0$. By the Serre duality and Corollary 3.1, we have $h^1(\mathcal{O}_X(-D))=0$.

We show that $h^1(\mathcal{O}_X(2)(-D))=0$. Let $\Delta$ be the fixed component of $|2C-D|$. We consider the case where $\Delta=\emptyset$. Since $(2C-D)^2=0$, $|2C-D|$ is base point free. By Remark 3.4, there exist a smooth irreducible curve $\tilde{D}$ on $X$ and a positive integer $m$ such that $m\tilde{D}\in|2C-D|$. Since $mC.\tilde{D}=C.(2C-D)=6$, if $m\geq2$, then $(m,C.\tilde{D})=(2,3)$, $(3,2)$ or $(6,1)$. Since $C.\tilde{D}=2P_a(\tilde{D})-2\equiv 0$ (mod 2), we have $m=3$. By the exact sequence 
$$0\longrightarrow\mathcal{O}_X\longrightarrow\mathcal{O}_X(\tilde{D})\longrightarrow\mathcal{O}_{\tilde{D}}(\tilde{D})\longrightarrow0,\leqno (4.1)$$
we have $h^0(\mathcal{O}_X(\tilde{D}))=2$. Since $|\tilde{D}-C|=\emptyset$ and $C.\tilde{D}=2$, this contradicts Remark 3.1, and hence, we have $m=1$. By the exact sequence (4.1), we have $h^0(\mathcal{O}_X(2)(-D))=2$. Since $|D-C|=\emptyset$ and $\chi(\mathcal{O}_X(2)(-D))=2$, we have $h^1(\mathcal{O}_X(2)(-D))=0$.

We consider the case where $\Delta\neq\emptyset$. We set $\tilde{D}=2C-D-\Delta$ and assume that $h^1(\mathcal{O}_X(2)(-D))\neq0$. Since $|D-C|=\emptyset$, we have $h^0(\mathcal{O}_X(2)(-D))\geq3$. By the ampleness of $C$, we have $C.\tilde{D}\leq5$. By Remark 3.1, we have $C.\tilde{D}=5$ and $h^0(\mathcal{O}_X(\tilde{D}))=3$. Since $C.\Delta=1$, by Lemma 3.1, $\Delta^2=-3$. Since $(2C-D)^2=0$, we have ${\tilde{D}}^2+2\tilde{D}.\Delta=3$. Since $\tilde{D}.\Delta\geq0$ and ${\tilde{D}}^2$ is a positive odd number, ${\tilde{D}}^2=1$ or 3. Since $C.(C-\tilde{D})=0$, we have $|C-\tilde{D}|=\emptyset$, and hence, $\chi(\mathcal{O}_X(\tilde{D}))\leq3$. This means that ${\tilde{D}}^2=1$. By Proposition 3.5, we can assume that $\tilde{D}$ is smooth and irreducible. By the exact sequence
$$0\longrightarrow\mathcal{O}_X\longrightarrow\mathcal{O}_X(\tilde{D})\longrightarrow\mathcal{O}_{\tilde{D}}(\tilde{D})\longrightarrow0,$$
we have $h^0(\mathcal{O}_{\tilde{D}}(\tilde{D}))=2$. However, since ${\tilde{D}}^2=1$ and $P_a(\tilde{D})>0$, this is a contradiction. Hence, $h^1(\mathcal{O}_X(2)(-D))=0$. Therefore, $D$ is an aCM curve. $\hfill\square$

$\;$

\noindent If a curve $D$ is linked to an aCM curve by a complete intersection of $X$ and a hypersurface in $\mathbb{P}^3$, $D$ is also aCM. Hence, we have the following assertion.

$\;$

\noindent{\bf{Corollary 4.1}} If $P_a(D)=4$ and $C.D=6$, then $D$ is an aCM curve.

\smallskip

\smallskip

{\it{Proof}}. Since $(D-C)^2=-7$ and $C.(D-C)=1$, by Lemma 3.1, we have $|D-C|=\emptyset$. Since $\chi(\mathcal{O}_X(2)(-D))=1$, we have $|2C-D|\neq\emptyset$. If we take $\tilde{D}\in|2C-D|$, $P_a(\tilde{D})=1$ and $C.\tilde{D}=4$. By Proposition 4.2, $\tilde{D}$ is an aCM curve. Hence, $D$ is also an aCM curve. $\hfill\square$

\begin{prop} If $P_a(D)=3$ and $C.D=5$, then $D$ is aCM. \end{prop}

{\it{Proof}}. Since $C.(D-C)=0$ and $D$ is not linearly equivalent to $C$, we have $|D-C|=\emptyset$ and $|C-D|=\emptyset$. Hence, $h^0(\mathcal{O}_X(D))\geq\chi(\mathcal{O}_X(D))=2$. Since $D^2=-1$, $|D|$ has a fixed component. We denote it by $\Delta$, and set $\tilde{D}=D-\Delta$. Then we have $C.\tilde{D}\leq4$. By Remark 3.1, we have $C.\tilde{D}=4$ and $h^0(\mathcal{O}_X(D))=2$, and hence, $h^1(\mathcal{O}_X(1)(-D))=0$. Since $|D-C|=\emptyset$, by Corollary 3.1, we have $h^1(\mathcal{O}_X(-D))=0$.

On the other hand, since $C.(2C-D)=5$, $(2C-D)^2=-1$, and $|D-C|=\emptyset$, we have $|2C-D|\neq\emptyset$. By the same reason as above, we have $h^1(\mathcal{O}_X(2)(-D))=0$ and $h^1(\mathcal{O}_X(3)(-D))=h^1(\mathcal{O}_X(D)(-2))=0$. Since $C.D+5<15$, by Proposition 4.1, $D$ is aCM. $\hfill\square$

$\;$

\noindent Non-aCM curves on $X$ can be explicitly written. From now on we give a necessary and sufficient condition for $D$ to be a non-aCM curve.

\begin{prop} Assume that $P_a(D)=6$ and $C.D=7$. Then the following conditions are equivalent.

\smallskip

\smallskip

\noindent{\rm{(a)}} $D$ is not an aCM curve.

\noindent{\rm{(b)}} There exist two lines $\Gamma_1$ and $\Gamma_2$ on $X$ with $\Gamma_1.\Gamma_2=0$ and $C+\Gamma_1+\Gamma_2\in|D|$. \end{prop}

{\it{Proof}}. (b) $\Longrightarrow$ (a). Since the member of $|D-C|$ is not 1-connected, $D$ is not an aCM curve.

(a) $\Longrightarrow$ (b). Since $(C-D).C=-2$, we have $|C-D|=\emptyset$. Since $(2C-D).C=3<4$, by Remark 3.1, if $|2C-D|\neq\emptyset$, we have $h^0(\mathcal{O}_X(2)(-D))=1$. In this case, by Theorem 1.1, $D$ is aCM. Hence, we have $|2C-D|=\emptyset$. By the Serre duality and the Riemann-Roch theorem, $h^0(\mathcal{O}_X(D)(-1))\geq\chi(\mathcal{O}_X(2)(-D))=1$. Since $(D-C).C=2$ and $(D-C)^2=-6$, by Lemma 3.2, we have the assertion. $\hfill\square$

$\;$

\noindent By Example 3.1, the non-aCM curve as in Proposition 4.4 can be constructed on the Fermat quintic hypersurface in $\mathbb{P}^3$.

\begin{prop} Assume that $P_a(D)=11$ and $C.D=10$. Then the following conditions are equivalent.

\smallskip

\smallskip

\noindent{\rm{(a)}} $D$ is not an aCM curve.

\noindent{\rm{(b)}} $|D-C|\neq\emptyset$ or $|3C-D|\neq\emptyset$. \end{prop}

\noindent{\bf{Remark 4.1}} In Proposition 4.5, since $(2C-D).C=0$ and $D$ is not linearly equivalent to $2C$, we have $|2C-D|=\emptyset$ and $|D-2C|=\emptyset$. Since $\chi(\mathcal{O}_X(2)(-D))=0$, the condition that $h^1(\mathcal{O}_X(2)(-D))=0$ is equivalent to $|D-C|=\emptyset$. Moreover, since $\chi(\mathcal{O}_X(3)(-D))=0$, $h^1(\mathcal{O}_X(3)(-D))=0$ if and only if $|3C-D|=\emptyset$.

$\;$

{\it{Proof of Proposition 4.5}}. (b) $\Longrightarrow$ (a). By Remark 4.1, the assertion is clear.

\noindent (a) $\Longrightarrow$ (b). Assume that $|D-C|=\emptyset$ and $|3C-D|=\emptyset$. By Remark 4.1, we have $h^1(\mathcal{O}_X(D)(-2))=h^1(\mathcal{O}_X(3)(-D))=0$. By the exact sequence
$$0\longrightarrow\mathcal{O}_X(D)(-2)\longrightarrow\mathcal{O}_X(D)(-1)\longrightarrow\mathcal{O}_C(D)(-1)\longrightarrow0,$$
we have $h^0(\mathcal{O}_C(D)(-1))=0$. By the exact sequence
$$0\longrightarrow\mathcal{O}_X(D)(-3)\longrightarrow\mathcal{O}_X(D)(-2)\longrightarrow\mathcal{O}_C(D)(-2)\longrightarrow0,$$
we have $h^1(\mathcal{O}_X(4)(-D))=h^1(\mathcal{O}_X(D)(-3))=0$. Since $C.(2C-D)=0$, if $|\mathcal{O}_C(2)(-D)|\neq\emptyset$, then we have $\mathcal{O}_C(D)\cong\mathcal{O}_C(2)$. This contradicts the fact that $h^0(\mathcal{O}_C(D)(-1))=0$. Hence, we have $h^0(\mathcal{O}_C(2)(-D))=0$.

On the other hand, since $|D-C|=\emptyset$, by Remark 4.1, $h^1(\mathcal{O}_X(2)(-D))=0$. By the exact sequence
$$0\longrightarrow\mathcal{O}_X(1)(-D)\longrightarrow\mathcal{O}_X(2)(-D)\longrightarrow\mathcal{O}_C(2)(-D)\longrightarrow0,$$
we have $h^1(\mathcal{O}_X(1)(-D))=0$. By Corollary 3.1, we have $h^1(\mathcal{O}_X(-D))=0$. Since $C.D+5<20$, by Proposition 4.1, $D$ is an aCM curve. $\hfill\square$

\begin{df} Let $\mathcal{F}$ be a coherent sheaf on $\mathbb{P}^3$. We say that $\mathcal{F}$ is $m$-{\rm{regular}}, if $h^i(\mathcal{F}(m-i))=0$ for each $i>0$. Moreover, we call the minimum number $m$ such that $\mathcal{F}$ is $m$-regular the regularity of $\mathcal{F}$.\end{df}

\noindent There exist non-aCM curves $D$ on $X$ satisfying the condition (b) as in Proposition 4.5. We prepare the following theorems to construct such curves lying on smooth hypersurfaces of degree $d\leq4$ in $\mathbb{P}^3$.

\begin{thm}{\rm{([3, p.100])}} If a coherent sheaf $\mathcal{F}$ on $\mathbb{P}^3$ is $m$-regular, then $\mathcal{F}(m)$ is globally generated. \end{thm}

\noindent The regularity of a curve in $\mathbb{P}^3$ is defined as the regularity of the ideal sheaf of it in $\mathbb{P}^3$. By the theorem of Bertini, we have the following assertion.

\begin{thm}{\rm{([2, Proposition 5.1 and Remark 5.2])}} If the regularity of a curve $D$ in $\mathbb{P}^3$ is $m$, then there exists a reduced and irreducible hypersurface $Y\in |\mathcal{I}_D(m)|$ which is smooth at any point of $Y\backslash\Supp(D)$. In particular, if $D$ is reduced and the embedding dimension at any point of $D$ is at most 2, then we can take such a surface $Y$ to be smooth.\end{thm}

\noindent In Theorem 4.2, the embedding dimension of $D$ at $P$ is defined as the dimension of the Zariski tangent space of $D$ at $P$. If $D$ is a reduced divisor on a smooth hypersurface in $\mathbb{P}^3$, then the embedding dimension at any point of $D$ is at most 2.

$\;$

\noindent{\bf{Example 4.1}} Let $Y$ be a smooth cubic surface in $\mathbb{P}^3$. Let $\pi:Y\longrightarrow\mathbb{P}^2$ be a blow up at six points $P_1,\cdots,P_6$ on $\mathbb{P}^2$ in general position. We set $E_i=\pi^{-1}(P_i)$ $(1\leq i\leq6)$, and let $l$ be the total transform of a line on $\mathbb{P}^2$ by $\pi$. Let $H_Y$ be the hyperplane class of $Y$. Then $-K_Y=H_Y$ and $H_Y\in|3l-\displaystyle\sum_{i=1}^6E_i|$. 

\noindent We take a reduced divisor $\tilde{D}\in|H_Y+E_1+E_2|$. Then $P_a(\tilde{D})=1$ and $H_Y.\tilde{D}=5$. Moreover, $\tilde{D}$ has regularity $\leq$ 5, in the sense of Castelnuovo-Mumford. Indeed, since $|H_Y-E_i|$ is an base point free pencil on $Y$ for $1\leq i\leq6$, the linear system $|4H_Y-E_1-E_2|$ is base point free and ample. Since $h^1(\mathcal{O}_Y(4)(-\tilde{D}))=h^1(K_Y(4)(-E_1-E_2))=0$, by the exact sequence
$$0\longrightarrow\mathcal{O}_{\mathbb{P}^3}(1)\longrightarrow\mathcal{I}_{\tilde{D}}(4)\longrightarrow\mathcal{O}_Y(4)(-\tilde{D})\longrightarrow0,$$
we have $h^1(\mathcal{I}_{\tilde{D}}(4))=0$. Since $h^2(\mathcal{O}_Y(3)(-\tilde{D}))=h^0(\mathcal{O}_Y(-3)(E_1+E_2))=0$, by the exact sequence
$$0\longrightarrow\mathcal{O}_{\mathbb{P}^3}\longrightarrow\mathcal{I}_{\tilde{D}}(3)\longrightarrow\mathcal{O}_Y(3)(-\tilde{D})\longrightarrow0,$$
we have $h^2(\mathcal{I}_{\tilde{D}}(3))=0$. Since $h^3(\mathcal{O}_{\mathbb{P}^3}(-1))=h^3(K_{\mathbb{P}^3}(3))=0$, by the similar reason as above, we have $h^3(\mathcal{I}_{\tilde{D}}(2))=0$. Hence, by Theorem 4.1, $\mathcal{I}_{\tilde{D}}(5)$ is generated by the global sections of it. Since $\tilde{D}$ has embedding dimension $\leq2$ at any point, by Theorem 4.2, there exists a smooth quintic hypersurface $X$ in $\mathbb{P}^3$ containing $\tilde{D}$. Let $D_1$ be a curve on $X$ which is linked to $\tilde{D}$ by the complete intersection $X\cap Y$. Then the degree of $D_1$ is 10 and $P_a(D_1)=11$. On the other hand, if we take a curve $D_2\in|\mathcal{O}_X(\tilde{D})(1)|$, $D_2$ also has the same degree and the arithmetic genus as $D_1$.

\begin{prop} Assume that $P_a(D)=3$ and $C.D=6$. Then the following conditions are equivalent.

\smallskip

\smallskip

\noindent{\rm{(a)}} $D$ is not an aCM curve.

\noindent{\rm{(b)}} One of the following cases occurs.

\smallskip

{\rm{(b$_1$)}} $|2C-D|\neq\emptyset$

{\rm{(b$_2$)}} There exist an effective divisor $\tilde{D}$ on $X$ with $P_a(\tilde{D})=3$ and $C.\tilde{D}=4$, and two lines $\Gamma_1$ and $\Gamma_2$ on $X$ such that $\Gamma_1.\Gamma_2=0$ and $\tilde{D}+\Gamma_1+\Gamma_2\in|D|$.

{\rm{(b$_3$)}} There exist an effective divisor $\tilde{D}$ on $X$ with $P_a(\tilde{D})=3$ and $C.\tilde{D}=4$, and a divisor $\Delta$ on $X$ such that $\Delta^2=-4$ and $\tilde{D}+\Delta\in|D|$. 

\end{prop}

\noindent{\bf{Remark 4.2}} Let $D$ be as in Proposition 4.6. Since $C.(D-C)=1$ and $(D-C)^2=-9$, by Lemma 3.1, $|D-C|=\emptyset$. Since $\chi(\mathcal{O}_X(2)(-D))=0$, by the Serre duality, $h^1(\mathcal{O}_X(2)(-D))=0$ if and only if $|2C-D|=\emptyset$. 

On the other hand, since $C.(C-D)=-1$, we have $|C-D|=\emptyset$. By the Serre duality, we have $h^2(\mathcal{O}_X(D))=0$. Since $\chi(\mathcal{O}_X(D))=1$, by the Riemann-Roch theorem, $h^1(\mathcal{O}_X(D))=0$ if and only if $h^0(\mathcal{O}_X(D))=1$.

$\;$

{\it{Proof of Proposition 4.6}}. (b) $\Longrightarrow$ (a). Assume $|2C-D|\neq\emptyset$. Then, by Remark 4.2, $D$ is not aCM. Assume that $D$ satisfies the condition (b$_2$) or (b$_3$) as in (b). An effective divisor $\tilde{D}$ on $X$ with $P_a(\tilde{D})=3$ and $C.\tilde{D}=4$ is a plane quartic. Indeed, since $C.(C-\tilde{D})=1$ and $(C-\tilde{D})^2=-3$, by Lemma 3.1, we have $|C-\tilde{D}|\neq\emptyset$. Since $|\tilde{D}|$ is a pencil on $X$, $h^0(\mathcal{O}_X(D))\geq2$. By Remark 4.2, this implies that $h^1(\mathcal{O}_X(D))\neq0$, and hence, $D$ is not aCM.

(a) $\Longrightarrow$ (b). Assume that $D$ does not satisfy the condition (b). Then we show that $D$ is aCM.  First of all, we show that $h^1(\mathcal{O}_X(D))=0$. Assume that $h^1(\mathcal{O}_X(D))\neq0$. By Remark 4.2, we have $h^0(\mathcal{O}_X(D))\geq2$. Since $D^2=-2$, $|D|$ has a fixed component $\Delta$. We set $\tilde{D}=D-\Delta$. Since $|D-C|=\emptyset$, we have $|\tilde{D}-C|=\emptyset$. By the ampleness of $C$, $C.\Delta\geq1$ and hence, $C.\tilde{D}\leq5$. Since $h^0(\mathcal{O}_X(\tilde{D}))\geq2$, by Remark 3.1, we have $C.\tilde{D}\geq4$. Therefore, we have
$$(C.\tilde{D},C.\Delta)=(5,1)\text{ or }(4,2).$$
We consider the case where $(C.\tilde{D},C.\Delta)=(5,1)$. By Lemma 3.1, we have $\Delta^2=-3$. Since $D^2=-2$ and $\tilde{D}$ is numerical effective, we have $(\tilde{D}^2,\tilde{D}.\Delta)=(1,0)$. Since $C.(C-\tilde{D})=0$ and $\tilde{D}$ is not linearly equivalent to $C$, we have $|C-\tilde{D}|=\emptyset$. By the Riemann-Roch theorem, we have $h^0(\mathcal{O}_X(\tilde{D}))\geq\chi(\mathcal{O}_X(\tilde{D}))=3$. Since $C.\tilde{D}=5$, by Remark 3.1, we have $h^0(\mathcal{O}_X(\tilde{D}))=3$ and $h^1(\mathcal{O}_X(\tilde{D}))=0$. By Proposition 3.5, we can assume that $\tilde{D}$ is smooth and irreducible. By the exact sequence
$$0\longrightarrow\mathcal{O}_X\longrightarrow\mathcal{O}_X(\tilde{D})\longrightarrow\mathcal{O}_{\tilde{D}}(\tilde{D})\longrightarrow0,$$
we have $h^0(\mathcal{O}_{\tilde{D}}(\tilde{D}))=2$. Since $\tilde{D}^2=1$ and $P_a(\tilde{D})>0$, this is a contradiction.

We consider the case where $(C.\tilde{D},C.\Delta)=(4,2)$. Since $h^0(\mathcal{O}_X(\tilde{D}))\geq2$, by Remark 3.1, we have $h^0(\mathcal{O}_X(\tilde{D}))=2$. Since $\tilde{D}^2\geq0$, we have $\chi(\mathcal{O}_X(\tilde{D}))\geq3$. Hence, by the Riemann-Roch theorem, we have $|C-\tilde{D}|\neq\emptyset$. Therefore, $\tilde{D}$ is a plane quartic and hence, we have $P_a(\tilde{D})=3$ and $C.\tilde{D}=4$. Since $C.\Delta=2$, $\Delta^2$ is even, and by Remark 3.2, we have $\Delta^2\leq-4$. By Lemma 3.2, we have $\Delta^2=-4,-6$, or $-12$, and if $\Delta^2=-12$, then $\Delta$ is a double line, that is, there exists a line $\Gamma$ on $X$ such that $\Delta\in|2\Gamma|$. Since $\tilde{D}^2=0$ and $D^2=-2$, we have the contradiction $2\tilde{D}.\Gamma=5$. If $\Delta^2=-4$, then $\Delta$ is a plane conic, and if $\Delta^2=-6$, then $\Delta$ is a divisor consisting of two skew lines on $X$. Hence, $D$ satisfies the condition (b$_2$) or (b$_3$) as in (b). This contradicts the first hypothesis. Therefore, we have $h^1(\mathcal{O}_X(D))=0$. Since $|D-C|=\emptyset$, by Corollary 3.1, we have $h^1(\mathcal{O}_X(-D))=0$. Moreover, since $|2C-D|=\emptyset$, by Remark 4.2, we have $h^1(\mathcal{O}_X(2)(-D))=0$. Since $C.D+5<15$, by Proposition 4.1, it is sufficient to show that $h^1(\mathcal{O}_X(3)(-D))=0$. First of all, since $|2C-D|=\emptyset$, by the exact sequence
$$0\longrightarrow\mathcal{O}_X(1)(-D)\longrightarrow\mathcal{O}_X(2)(-D)\longrightarrow\mathcal{O}_C(2)(-D)\longrightarrow0,$$
we have $h^0(\mathcal{O}_C(2)(-D))=0$. On the other hand, $h^0(\mathcal{O}_C(D)(-1))=0$. In fact, since $C.(D-C)=1$, if $h^0(\mathcal{O}_C(D)(-1))>0$, then there exists a point $P\in C$ such that $\mathcal{O}_C(D)\cong\mathcal{O}_C(1)(P)$. This means that $\mathcal{O}_C(2)(-D)\cong\mathcal{O}_C(1)(-P)$. However, this is a contradiction. Since $h^1(\mathcal{O}_X(D)(-1))=h^1(\mathcal{O}_X(2)(-D))=0$, by the exact sequence
$$0\longrightarrow\mathcal{O}_X(D)(-2)\longrightarrow\mathcal{O}_X(D)(-1)\longrightarrow\mathcal{O}_C(D)(-1)\longrightarrow0,$$
we have $h^1(\mathcal{O}_X(3)(-D))=h^1(\mathcal{O}_X(D)(-2))=0$. Therefore, $D$ is aCM. $\hfill\square$

$\;$

\noindent First of all, we construct a non-aCM curve $D$ satisfying the condition (b$_1$)  as in Proposition 4.6 (b), by using a curve lying on a smooth quadric in $\mathbb{P}^3$. 

$\;$

\noindent{\bf{Example 4.2}} Let $Z$ be a smooth quadric in $\mathbb{P}^3$. Note that since $Z\cong \mathbb{P}^1\times\mathbb{P}^1$, if we let $L_1$ and $L_2$ be the classes of two skew lines on $Z$, then $\Pic(Z)=\mathbb{Z}L_1\oplus\mathbb{Z}L_2$. Let $H_Z$ be the hyperplane class of $Z$. Then $H_Z=L_1+L_2$, and $K_Z=-2H_Z$. Since $|H_Z+2L_2|$ is base point free and big, by the theorem of Bertini, we can take a smooth irreducible curve $\tilde{D}\in |H_Z+2L_2|$. Then $P_a(\tilde{D})=0$ and $H_Z.\tilde{D}=4$. Moreover, $\tilde{D}$ has regularity $\leq5$. Indeed, since $|5H_Z-2L_2|$ is base point free and big, we have $h^1(\mathcal{O}_Z(4)(-\tilde{D}))=h^1(K_Z(5)(-2L_2))=0$. By the exact sequence
$$0\longrightarrow\mathcal{O}_{\mathbb{P}^3}(2)\longrightarrow\mathcal{I}_{\tilde{D}}(4)\longrightarrow\mathcal{O}_Z(4)(-\tilde{D})\longrightarrow0,$$
we have $h^1(\mathcal{I}_{\tilde{D}}(4))=0$. Since $h^2(\mathcal{O}_Z(3)(-\tilde{D}))=h^0(\mathcal{O}_Z(-5)(\tilde{D}))=0$, by the exact sequence
$$0\longrightarrow\mathcal{O}_{\mathbb{P}^3}(1)\longrightarrow\mathcal{I}_{\tilde{D}}(3)\longrightarrow\mathcal{O}_Z(3)(-\tilde{D})\longrightarrow0,$$
we have $h^2(\mathcal{I}_{\tilde{D}}(3))=0$. Moreover, since $h^3(\mathcal{O}_{\mathbb{P}^3})=0$, we have $h^3(\mathcal{I}_{\tilde{D}}(2))=0$. Hence, by Theorem 4.2, there exists a smooth quintic hypersurface $X$ in $\mathbb{P}^3$ containing $\tilde{D}$. If we let $D$ be a non-zero effective divisor on $X$ which is linked to $\tilde{D}$ by the complete intersection $X\cap Z$, we have $P_a(D)=3$ and $C.D=6$.

$\;$

\noindent On the other hand, we can give an example of a non-aCM curve $D$ satisfying the condition (b$_2$) or (b$_3$) as in Proposition 4.6 (b), by using two skew lines on the Fermat quintic hypersurface in $\mathbb{P}^3$.

$\;$

\noindent {\bf{Example 4.3}} Let $X$ be the quintic hypersurface in $\mathbb{P}^3$ defined by the equation $x_0^5+x_1^5+x_2^5+x_3^5=0$. Let $L_1$ and $L_2$ be lines on $X$ defined by the equations $x_0+x_1=x_2+x_3=0$ and $x_0+x_2=x_1+\xi x_3=0$, respectively. Here, $\xi$ is a primitive fifth root of unity. Note that $L_1$ and $L_2$ do not intersect (see Example 3.1). Moreover, let $\tilde{C}$ be the hyperplane section of $X$ defined by $x_1+x_2=0$, and let $\Gamma$ be the line defined by the equation $x_1+x_2=x_0+\xi^4x_3=0$. Since $L_i.\Gamma=0$ and $\tilde{C}.L_i=1$ ($i=1,2$), $\tilde{C}-\Gamma$ is a plane quartic curve such that $(\tilde{C}-\Gamma).L_i=1$ ($i=1,2$). Hence, if we set $\tilde{D}=\tilde{C}-\Gamma$ and  $D=\tilde{D}+L_1+L_2$, then $D$ is a non-aCM curve corresponding to the case (b$_2$) as in Proposition 4.6 (b).

On the other hand, let $L_3$ be the line on $X$ defined by the equation $x_0+x_1=x_2+\xi^4x_3=0$. Since $L_3.\Gamma=1$, we have $L_3.\tilde{D}=0$. The divisor $L_1+L_3$ is a plane conic which is contained in the hyperplane section of $X$ defined by the equation $x_0+x_1=0$. Since $(L_1+L_3).\tilde{D}=1$, if we set $D=\tilde{D}+L_1+L_3$, $D$ is a non-aCM curve corresponding to the case (b$_3$) as in Proposition 4.6 (b).

$\;$

\noindent{\bf{Corollary 4.2}} If $P_a(D)=C.D=9$, then the following conditions are equivalent.

\smallskip

\smallskip

\noindent{\rm{(a)}} $D$ is not an aCM curve.

\noindent{\rm{(b)}} One of the following cases occurs.

\smallskip

{\rm{(b$_1$)}} $|D-C|\neq\emptyset$

{\rm{(b$_2$)}} There exist an effective divisor $\tilde{D}$ on $X$ with $P_a(\tilde{D})=3$ and $C.\tilde{D}=4$, and two lines $\Gamma_1$ and $\Gamma_2$ on $X$ such that $\Gamma_1.\Gamma_2=0$ and $\tilde{D}+\Gamma_1+\Gamma_2\in|3C-D|$.

{\rm{(b$_3$)}} There exist an effective divisor $\tilde{D}$ on $X$ with $P_a(\tilde{D})=3$ and $C.\tilde{D}=4$, and a divisor $\Delta$ on $X$ such that $\Delta^2=-4$ and $\tilde{D}+\Delta\in|3C-D|$. 

$\;$

{\it{Proof}}. First of all, since $C.(D-2C)=-1$, we have $|D-2C|=\emptyset$. Since $\chi(\mathcal{O}_X(3)(-D))=1$, we have $|3C-D|\neq\emptyset$. If we take $\tilde{D}\in |3C-D|$, since $P_a(\tilde{D})=3$ and $C.\tilde{D}=6$, by Proposition 4.6, we have the assertion. $\hfill\square$

$\;$

\noindent Note that we can construct non-aCM curves satisfying the condition (b) as in Corollary 4.2, by Example 4.2, Example 4.3, and the proof of Corollary 4.2.

\begin{prop} Assume that $P_a(D)=5$ and $C.D=7$. Then the following conditions are equivalent.

\smallskip

\smallskip

\noindent{\rm{(a)}} $D$ is not an aCM curve.

\noindent{\rm{(b)}} There exist effective divisors $\Gamma_1$ and $\Gamma_2$ on $X$ with $C.\Gamma_i=i$ {\rm{(}}$i=1,2${\rm{)}}, $\Gamma_1^2=-3$, $\Gamma_2^2=-4$, $\Gamma_1.\Gamma_2=0$, and $\Gamma_1+\Gamma_2\in|2C-D|$. \end{prop}

{\it{Proof}}. (b) $\Longrightarrow$ (a). Since the member of $|2C-D|$ is not 1-connected, $D$ is not aCM.

(a) $\Longrightarrow$ (b). Since $(D-C)^2=-8$ and $C.(D-C)=2$, by Lemma 3.2, we have $|D-C|=\emptyset$. Assume that $D$ does not satisfy the condition (b). Then we show that $D$ is aCM. First of all, we show that $h^1(\mathcal{O}_X(D))=0$. Since $C.(C-D)=-2$, we have $|C-D|=\emptyset$. Assume that $h^1(\mathcal{O}_X(D))\neq0$. Then, by the Riemann-Roch theorem, we have $h^0(\mathcal{O}_X(D))\geq3$. Let $\Delta$ be the fixed component of $|D|$, and assume that it is not empty. We set $\tilde{D}=D-\Delta$. Since $C.\Delta\geq1$, we have $C.\tilde{D}\leq6$. By Remark 3.1, we have $C.\tilde{D}-2\geq h^0(\mathcal{O}_X(\tilde{D}))\geq3$. Hence, $C.\tilde{D}=5$ or 6 (resp. $C.\Delta=2$ or 1). Since $C.(C-\tilde{D})\leq0$ and $\tilde{D}$ is not linearly equivalent to $C$, we have $|C-\tilde{D}|=\emptyset$. Therefore, since $C.\tilde{D}-2\geq\chi(\mathcal{O}_X(\tilde{D}))$, we have
$$3C.\tilde{D}-14\geq\tilde{D}^2.\leqno(4.2)$$

We consider the case where $C.\tilde{D}=6$. Since $C.\Delta=1$, by Lemma 3.1, we have $\Delta^2=-3$. Moreover, $\tilde{D}^2$ is even and, by the inequality (4.2), we have $0\leq\tilde{D}^2\leq4$. Assume that $\tilde{D}^2=4$. Since $C.(\tilde{D}-C)=1$ and $(\tilde{D}-C)^2=-3$, by Lemma 3.1, $|\tilde{D}-C|\neq\emptyset$. This contradicts the fact that $|D-C|=\emptyset$. 

Assume that $\tilde{D}^2=2$. Since $|\tilde{D}-C|=\emptyset$, by the Riemann-Roch theorem, we have $h^0(\mathcal{O}_X(2)(-\tilde{D}))\geq\chi(\mathcal{O}_X(2)(-\tilde{D}))=2$. Hence, we take $D_0\in |2C-\tilde{D}|$. Since $|D_0-C|=|C-\tilde{D}|=\emptyset$, by Remark 3.1, we have $C.D_0=4$ and $h^0(\mathcal{O}_C(D_0))=2$. Hence, $|\mathcal{O}_C(D_0)|$ is a gonality pencil on $C$. On the other hand, since $D_0^2<0$, $|D_0|$ has a fixed component. However, since $C$ is ample, this is a contradiction. 

Assume that $\tilde{D}^2=0$. Since $|\tilde{D}|$ is base point free, by Remark 3.4, there exist a smooth irreducible curve $D_0$ and a positive integer $m$ such that $mD_0\in|\tilde{D}|$. Since $C.D_0=2P_a(D_0)-2\equiv0$ (mod 2), $m=1$ or 3. If $m=1$, then by the exact sequence
$$0\longrightarrow\mathcal{O}_X\longrightarrow\mathcal{O}_X(D_0)\longrightarrow\mathcal{O}_{D_0}(D_0)\longrightarrow0,\leqno (4.3)$$
we have $h^0(\mathcal{O}_X(\tilde{D}))=h^0(\mathcal{O}_X(D_0))=2$. Since $h^0(\mathcal{O}_X(D))\geq3$, this is a contradiction. If $m=3$, we have $C.D_0=2$.  By the exact sequence (4.3), $h^0(\mathcal{O}_X(D_0))=2$. However, by Remark 3.1, this is a contradiction. 

We consider the case where $C.\tilde{D}=5$. By the inequality (4.2), we have $\tilde{D}^2\leq1$. Since $\tilde{D}^2+5=2P_a(\tilde{D})-2$, $\tilde{D}^2$ is an odd number. Hence, we have $\tilde{D}^2=1$. By Proposition 3.5, we can assume that $\tilde{D}$ is smooth and irreducible. Since $P_a(\tilde{D})>0$, we have $h^0(\mathcal{O}_{\tilde{D}}(\tilde{D}))=1$. Hence, by the exact sequence
$$0\longrightarrow\mathcal{O}_X\longrightarrow\mathcal{O}_X(\tilde{D})\longrightarrow\mathcal{O}_{\tilde{D}}(\tilde{D})\longrightarrow0,\leqno (4.4)$$
we have $h^0(\mathcal{O}_X(\tilde{D}))=2$. However, this contradicts the fact that $h^0(\mathcal{O}_X(D))\geq3$. By the above argument, we have $\Delta=\emptyset$. Since $D^2=1$, by Proposition 3.5, we can assume that $D$ is smooth and irreducible, and hence, by the similar reason as above, we have a contradiction. Therefore, we have $h^1(\mathcal{O}_X(D))=0$. Since $|D-C|=\emptyset$, by Corollary 3.1, we have $h^1(\mathcal{O}_X(-D))=0$. 

Next, we show that $h^1(\mathcal{O}_X(2)(-D))=0$. Assume that $h^1(\mathcal{O}_X(2)(-D))\neq0$. Since $\chi(\mathcal{O}_X(2)(-D))=0$ and $|D-C|=\emptyset$, we have $|2C-D|\neq\emptyset$. We take $\tilde{D}\in|2C-D|$. Since $P_a(\tilde{D})=-1$, there exists a non-trivial effective decomposition $\tilde{D}=\Gamma_1+\Gamma_2$ such that $\Gamma_1.\Gamma_2\leq0$. Since $C.\tilde{D}=3$, we may assume that $C.\Gamma_i=i$ ($i=1,2$). Then we note that $\Gamma_1^2=-3$, by Lemma 3.1. Hence, we have $\Gamma_2^2\geq\Gamma_2^2+2\Gamma_1.\Gamma_2=-4$. By Remark 3.2, we have $\Gamma_2^2=-4$. This means that $D$ satisfies the condition (b) and hence, we have $h^1(\mathcal{O}_X(2)(-D))=0$.

Since $C.D+5<15$, it is sufficient to show that $h^1(\mathcal{O}_X(3)(-D))=0$. Since $h^1(\mathcal{O}_X(D))=0$, by the Riemann-Roch theorem, we have $h^0(\mathcal{O}_X(D))=2$. Since $|D-C|=\emptyset$ and $h^1(\mathcal{O}_X(D)(-1))=h^1(\mathcal{O}_X(2)(-D))=0$, by the exact sequence
$$0\longrightarrow\mathcal{O}_X(D)(-1)\longrightarrow\mathcal{O}_X(D)\longrightarrow\mathcal{O}_C(D)\longrightarrow0,$$
we have $h^0(\mathcal{O}_C(D))=2$. Since $h^0(\mathcal{O}_C(1))=3$, we have $h^0(\mathcal{O}_C(D)(-1))=0$. Moreover, by the exact sequence
$$0\longrightarrow\mathcal{O}_X(D)(-2)\longrightarrow\mathcal{O}_X(D)(-1)\longrightarrow\mathcal{O}_C(D)(-1)\longrightarrow0,$$
we have $h^1(\mathcal{O}_X(3)(-D))=h^1(\mathcal{O}_X(D)(-2))=0$. Hence, $D$ is aCM. $\hfill\square$

$\;$

\noindent A non-aCM curve satisfying the condition (b) as in Proposition 4.7 also can be constructed as a divisor on the Fermat quintic hypersurface in $\mathbb{P}^3$.

$\;$

\noindent {\bf{Example 4.4}} Let $X$ be the quintic hypersurface in $\mathbb{P}^3$ defined by the equation $x_0^5+x_1^5+x_2^5+x_3^5=0$. Let $C_1$ and $C_2$ be the hyperplane sections defined by the equations $x_0+x_1=0$ and $x_0+x_2=0$, respectively. Moreover, let $\Gamma_1$ be the line on $X$ defined by the equation $x_0+x_1=x_2+x_3=0$, and let $\Gamma_2$ be the conic on $X$ defined by the equation $x_0+x_2=(x_1+\xi x_3)(x_1+\xi^2x_3)=0$, where $\xi$ is a primitive fifth root of unity. Since any irreducible component of $\Gamma_2$ does not intersect $\Gamma_1$, we have $\Gamma_1.\Gamma_2=0$. Hence, if we set $D=C_1+C_2-\Gamma_1-\Gamma_2$, then $D$ is a non-aCM curve corresponding to the case (b) as in Proposition 4.7.

$\;$

\noindent{\bf{Corollary 4.3}} Assume that $P_a(D)=7$ and $C.D=8$. Then the following conditions are equivalent.

\smallskip

\smallskip

\noindent{\rm{(a)}} $D$ is not an aCM curve.

\noindent{\rm{(b)}} There exist effective divisors $\Gamma_1$ and $\Gamma_2$ on $X$ with $C.\Gamma_i=i$ {\rm{(}}$i=1,2${\rm{)}}, $\Gamma_1^2=-3$, $\Gamma_2^2=-4$, $\Gamma_1.\Gamma_2=0$, and $\Gamma_1+\Gamma_2\in|D-C|$.

$\;$

{\it{Proof}}. First of all, since $C.(D-2C)=-2$, we have $|D-2C|=\emptyset$. Since $\chi(\mathcal{O}_X(3)(-D))=2$, we have $|3C-D|\neq\emptyset$. If we take $\tilde{D}\in |3C-D|$, since $P_a(\tilde{D})=5$ and $C.\tilde{D}=7$, by Proposition 4.7, we have the assertion. $\hfill\square$

$\;$

\noindent If $D$ is the divisor on the Fermat quintic hypersurface $X$ in $\mathbb{P}^3$ as in Example 4.4, we can construct a non-aCM curve satisfying the condition (b) as in Corollary 4.3 as a curve which is linked to $D$ by the complete intersection of $X$ and a cubic hypersurface in $\mathbb{P}^3$ containing $D$.

\begin{prop} Assume that $P_a(D)=2$ and $C.D=5$. Then the following conditions are equivalent.

\smallskip

\smallskip

\noindent {\rm{(a)}} $D$ is not an aCM curve.

\smallskip

\smallskip

\noindent {\rm{(b)}} There exist a line $\Gamma$ and an effective divisor $\tilde{D}$ on $X$ with $C.\tilde{D}=4$, $P_a(\tilde{D})=3$, $\tilde{D}.\Gamma=0$, and $\tilde{D}+\Gamma\in|D|$ or $|2C-D|$. \end{prop}

{\it{Proof}}. (b) $\Longrightarrow$ (a). If $D$ satisfies the condition (b), the members of $|D|$ or $|2C-D|$ are not 1-connected. Hence, $D$ is not aCM.

(a) $\Longrightarrow$ (b). Since $(C-D).C=0$ and $P_a(D)=2$, we have $|C-D|=|D-C|=\emptyset$. Assume that $D$ does not satisfy the condition (b). Then we show that $D$ is an aCM curve. First of all, we show that $h^1(\mathcal{O}_X(D))=0$. Assume that $h^1(\mathcal{O}_X(D))\neq0$. Since $\chi(\mathcal{O}_X(D))=1$, we have $h^0(\mathcal{O}_X(D))\geq2$. Since $D^2=-3$, the fixed component of $|D|$ is not empty. Hence, let $\Gamma$ be the fixed component of $|D|$, and we take $\tilde{D}\in|D-\Gamma|$. Since $C.\tilde{D}\leq4$, by Remark 3.1, we have $C.\tilde{D}=4$ and $C.\Gamma=1$. Hence, $\Gamma$ is a line on $X$. Since $\tilde{D}^2+2\tilde{D}.\Gamma=0$, $\tilde{D}^2\geq0$, and $\tilde{D}.\Gamma\geq0$, we have $\tilde{D}^2=\tilde{D}.\Gamma=0$. In this case, $D$ satisfies the condition (b). Hence, we have $h^1(\mathcal{O}_X(D))=0$. Since $|D-C|=\emptyset$, by Corollary 3.1, we have $h^1(\mathcal{O}_X(-D))=0$. Since $\chi(\mathcal{O}_X(2)(-D))=1$, by the Serre duality and the Riemann-Roch theorem, we have $|2C-D|\neq\emptyset$. If we take a member $D_0\in |2C-D|$, then $P_a(D_0)=2$ and $C.D_0=5$. Hence, by symmetry, we have $h^1(\mathcal{O}_X(2)(-D))=0$ and $h^1(\mathcal{O}_X(3)(-D))=h^1(\mathcal{O}_X(D)(-2))=0$. Since $C.D+5<15$, by Proposition 4.1, $D$ is aCM. $\hfill\square$

$\;$

\noindent {\bf{Example 4.5}} Let $X$ be the quintic hypersurface in $\mathbb{P}^3$ defined by the equation $x_0^5+x_1^5+x_2^5+x_3^5=0$. Let $\Gamma$ and $\tilde{\Gamma}$ be lines on $X$ defined by the equations $x_0+x_1=x_2+x_3=0$ and $x_0+x_2=x_1+x_3=0$, respectively. Note that $\Gamma$ and $\tilde{\Gamma}$ intersect at one point. If we let $\tilde{C}$ be the hyperplane section of $X$ defined by the equation $x_0+x_2=0$, then $\tilde{C}.\Gamma=1$. Hence, we have $(\tilde{C}-\tilde{\Gamma}).\Gamma=0$. We set $\tilde{D}=\tilde{C}-\tilde{\Gamma}$ and $D=\tilde{D}+\Gamma$. Then $P_a(D)=2$ and $\tilde{C}.D=5$. Let $C$ be the hyperplane section defined by the equation $x_0+x_1=0$, and we set $D_0=C+\tilde{C}-D$. Then $D_0$ has the same arithmetic genus and degree as $D$.

$\;$

\noindent The non-aCM curves satisfying the condition (b) as in Proposition 4.8 (b) can not be constructed on any smooth quartic hypersurface in $\mathbb{P}^3$. In fact, a curve $\tilde{D}$ with $P_a(\tilde{D})=3$ and $\deg \tilde{D}=4$ on a smooth quartic hypersurface $Y$ in $\mathbb{P}^3$ is a hyperplane section of $Y$, and hence, any curve on $Y$ intersects $\tilde{D}$. Proposition 4.8 means that a curve of arithmetic genus 2 and degree 5 is not necessarily aCM if it does not lie on any smooth hypersurface of degree $d\leq4$ in $\mathbb{P}^3$. 

$\;$

\noindent {\bf{Acknowledgements}}

\smallskip

\smallskip

\noindent I would like to thank Akira Ohbuchi and Jiryo Komeda for the opportunity to talk about this topic at the 19 th symposium on algebraic curves theory. The author also would like to thank the referee for some helpful comments.

\end{document}